\documentstyle[11pt]{article}
\begin{document}
\newtheorem{Def}{Definition}[section]
\newtheorem{thm}{Theorem}[section]
\newtheorem{lem}{Lemma}[section]
\newtheorem{rem}{Remark}[section]
\newtheorem{prop}{Proposition}[section]
\newtheorem{cor}{Corollary}[section]
\def\Balpha{  {\bar \alpha} }
\def\Aalpha{  {\mu} }
\def\av{{\int \hspace{-2.25ex}-} }
\title
{On the $C^1$ regularity
of solutions to  divergence form elliptic systems
with Dini-continuous coefficients
}
\author{YanYan Li\thanks{Partially
 supported by
an  NSF grant.}\\
Department of Mathematics\\
Rutgers University\\
110 Frelinghuysen Rd.\\
Piscataway, NJ 08854-8019}

\date{}
\maketitle

\bigskip

\bigskip

\centerline{ Dedicated to Haim Brezis with friendship and admiration}

\bigskip

\input { amssym.def}

\paragraph{Abstract}
{We prove $C^1$ regularity of solutions to  divergence form elliptic systems
with Dini-continuous coefficients.}

\bigskip

\bigskip

\setcounter{section}{1}

This note  addresses a question raised to
the author  by Haim Brezis, in connection with his 
solution of a conjecture of Serrin
 concerning
 divergence form second order
elliptic equations (see \cite{B} and \cite{B1}).
If the coefficients of the equations (or systems) are
H\"older
continuous, then  $H^1$ solutions are known to have
H\"older
continuous first derivatives.
The question is what minimal regularity assumption
 of the
coefficients
 would guarantee $C^1$ regularity of
all  $H^1$ solutions.  See \cite{JMV} for answers to some other
related questions of Haim.

 Consider 
the elliptic system for vector-valued functions
$u=(u^1, \cdots, u^N)$,
$$
\partial_\alpha( A_{ij}^{\alpha\beta}(x)\partial_\beta u^j)
=0,  \quad  \mbox{in}\ B_4, \qquad i=1, 2, \cdots, N,
$$
where
 $B_4$ is the ball  in  
  $\Bbb R^n$  of  $4$ centered at the origin.
The coefficients
$\{A_{ij}^{\alpha\beta}\}$ satisfy, for some positive constants
  $\Lambda$ and  $\lambda$,
 \begin{equation}
|A_{ij}^{\alpha\beta}(x)|\le  \Lambda,
\quad x\in B_4,
\label{3}
\end{equation}
\begin{equation}
\int_{B_4}  A^{\alpha\beta}_{ij}(x)\partial_\alpha
\eta ^i\partial_\beta\eta^j\ge
\lambda \int _{B_4} |\nabla \eta|^2,
\qquad\forall\ \eta\in H^1_0(B_4, \Bbb R^N),
\label{4}
\end{equation}
and
\begin{equation}
\int_0^1 r^{-1}\bar \varphi(r)dr<\infty,
\label{z1}
\end{equation}
where
\begin{equation}
\bar \varphi(r):=
\sup_{x\in B_3}
( \av_{B_r(x)} |A- A(x)|^2)^{\frac 12}.
\label{pp}
\end{equation}

\medskip

\noindent{\bf Main Theorem}.\
{\it Suppose that $\{A_{ij}^{\alpha\beta}\}$
satisfy the above assumptions, and
$u\in H^1(B_4, \Bbb R^N)$ is a solution of the 
elliptic system.  Then
$u$ is  $C^1$ in $B_1$.
}

\medskip

\noindent{\bf Remark.}\ For elliptic equations with
coefficients
satisfying $\alpha-$increasing  Dini  conditions,
a proof of the $C^1$ regularity of $u$ can be found,
in e.g. 
[\cite{Lie}, Theorem 5.1] as pointed out in \cite{B} and \cite{B1}. 

\medskip

\noindent{\bf Question.}\
{\it If we replace $\bar\varphi$ in   (\ref{z1}) by}
\begin{equation}
\hat \varphi(r):=
\sup_{x\in B_3}
( \av_{B_r(x)} |A- \overline A_{ B_r(x) }|^2)^{\frac 12},
\label{ppp}
\end{equation}
with $\displaystyle{   \overline A_{ B_r(x) }
:=  \av_{B_r(x)}  A  }$, 
does the conclusion of the Main Theorem still hold?

\bigskip


Let $B_r(x) \subset \Bbb R^n$ denote the
ball of radius $r$ and centered at $x$. We often write 
$B_r$  for $B_r(0)$, and $rB_1$ for  $B_r$.
Consider 
elliptic systems \begin{equation}
\partial_\alpha( A_{ij}^{\alpha\beta}(x)\partial_\beta u^j)
=h_i+\partial_\beta f^\beta_i, \quad  \mbox{in}\ B_4,
\qquad i=1,\cdots,N,
\label{2}
\end{equation}
where 
 $\alpha,\beta$ are summed from $1$ to $n$, while $i,j$ are
summed from $1$ to $N$.
The coefficients 
$\{A_{ij}^{\alpha\beta}\}$, often denoted
by $A$, 
 satisfy, for some positive constants $\Lambda$ and $\lambda$,
(\ref{3}), (\ref{4}) and (\ref{z1}), with
$\bar\varphi$ given by (\ref{pp}).
\begin{thm}
For $B_4\subset \Bbb R^n$,
$n\ge 1$,    let
  $A$, 
$\Lambda$, $\lambda$, 
$\varphi$  be as above,  $\{h_i\}, 
\{f^\beta_i\}\in C^\alpha(B_4)$
for some $\alpha>0$,  and let $u\in H^1(B_4, \Bbb R^N)$,
$N\ge 1$, be a solution of (\ref{2}).
Then $u\in C^1(B_1)$.  Moreover, the modulus of continuity of
$\nabla u$ in $B_1$ can be controlled in terms of
$\bar \varphi$, $n$, $N$, $\Lambda$, $\lambda$, $\alpha$, 
$\|h\|_{  C^\alpha(B_2) }$ and $[f]_{ C^\alpha(B_2)}$.
\label{thm1}
\end{thm}
\begin{rem} 
Assumption (\ref{z1}) is weaker than
$A$ being Dini-continuous.
\end{rem}
\begin{rem} The conclusion of Theorem \ref{thm1} 
still holds (the dependence
on $\alpha$, $\|h\|_{  C^\alpha(B_2) }$ and $[f]_{ C^\alpha(B_2)}$
is changed accordingly) if $\{h_i\}\in
L^p(B_4)$ for some $p>n$, and 
$f$ satisfies
$$
\int_0^1 r^{-1} \bar\psi(r)dr<\infty,
\ \ \mbox{where}\ 
\bar\psi(r):=
\sup_{x\in B_3}
 ( \av_{B_r(x)} |f- f(x)|^2)^{\frac 12}.
$$
\end{rem}

\begin{rem}  This note was written
in 2008.  It was intended to be published after 
having an answer to the question raised above.
\end{rem} 

Theorem \ref{thm1} follows from the following two propositions.

\begin{prop}
For $B_4\subset \Bbb R^n$,
$n\ge 1$,
 let $\Lambda$, $\lambda$, $N$  be as above,
 and let
 $A$ satisfy
(\ref{3}), (\ref{4}), and
\begin{equation}
( \av_{B_r} |A- A(0)|^2)^{\frac 12}
\le \varphi(r),\quad 0<r<1,
\label{dini}
\end{equation}
for some  non-negative function  $\varphi$  on $(0,1)$ satisfying,
for some $\mu>1$,
\begin{equation}
\max_{ r/2 \le s\le 2r} \varphi(s)\le \mu \varphi(r),
\qquad \int_0^1 r^{-1}\varphi(r)dr<\infty.
\label{BBB}
\end{equation}
Assume that $h,
f\in  C^\alpha(B_4)$ for some $\alpha>0$,  and 
 $u\in H^1(B_4, \Bbb R^N)$
 is  a solution of (\ref{2}).
Then there exist $a\in \Bbb R$ and
$b\in \Bbb R^n$ such that
\begin{equation}
\av_{B_r}| u(x)-[a+b\cdot x]|dx
\le  r \delta(r)[ \|u\|_{ L^2(B_2) }
+ \|h\|_{ C^\alpha(B_2) }+ [f]_{ C^\alpha(B_2) }],\ \ \ 
\forall\ 0<r<1,
\label{ee1}
\end{equation}
where $ 
\delta(r)$, depending only on $
\varphi$, $n, \lambda, \Lambda, N, \mu,  \alpha$,
satisfies $\lim _{r\to 0}\delta(r)=0$.
\label{thmmain}
\end{prop}

\begin{prop} Let $u$ be  a Lebesgue integrable
 function on $B_1\subset \Bbb R^n$, $n\ge 1$,
and let 
$\delta(r)$ be a monotonically
increasing positive function defined on $(0, 1)$ satisfying
$
\displaystyle{\lim_{r\to 0}\delta(r)=0.}$
Assume that
  for every  $\bar x\in B_{1/4}$, there exist
$a(\bar x)\in \Bbb R$, $b(\bar x)\in \Bbb R^n$ such that
\begin{equation}
\av _{ B_r(\bar x)  }|u(x)- [a(\bar x)+b(\bar x)
\cdot (x -\bar x)]|
dx
\le r\delta(r), 
 \qquad \forall\ 0<r<1/2.
\label{abab3}
\end{equation}
Then
$u$, after changing its
values on a zero Lebesgue measure set, belongs to
$C^1(B_{1/4})$, with $u\equiv a$
and $\nabla u\equiv b$.  Moreover, for 
 some dimensional constant $C$,
\begin{equation}
|\nabla u(x)-\nabla u(y)|\le C \delta(4|x-y|),\ \ \ 
\forall\ x, y\in B_{1/4}.
\label{abab4}
\end{equation}
\label{lem1}
\end{prop}

Similar results hold for Dirichlet problem: 
Let $\Omega\subset 
\Bbb R^n$, $n\ge 1$, be a domain with smooth boundary,
let $\Lambda$ and $\lambda$ be positive constants, 
and let
$A$ satisfy, for $N\ge 1$,  
$$
|A_{ij}^{\alpha\beta}(x)|\le  \Lambda,\ x\in \Omega,
$$
$$
\int_\Omega  A^{\alpha\beta}_{ij}(x)\partial_\alpha
\eta ^i\partial_\beta\eta^j\ge
\lambda \int _\Omega |\nabla \eta|^2,
\qquad\forall\ \eta\in H^1_0(\Omega,  \Bbb R^N),
$$ 
$$
\int_0^1 r^{-1}\bar \psi(r)dr<\infty,
$$
where
$$
\bar\psi(r):=
\sup_{x\in \Omega}
( \av_{B_r(x)\cap \Omega} |A- A(x)|^2)^{\frac 12}.
$$
Consider
\begin{eqnarray*}
\partial_\alpha( A_{ij}^{\alpha\beta}(x)\partial_\beta u^j)
&=&h_i+\partial_\beta f^\beta_i, \quad  \mbox{in}\ \Omega,
\qquad i=1,\cdots,N,\\
u&=&g, \qquad\mbox{on}\ \partial \Omega.
\end{eqnarray*}
\begin{thm}\ (\cite{L})\  Assume the above, and  let $h, f\in 
C^\alpha(\overline \Omega)$ and $g\in C^{1, \alpha}(\partial \Omega)$
for some $\alpha>0$.
  Then an $H^1(\Omega, \Bbb R^N)$ solution
$u$ to the above Dirichlet problem is in $C^1(\overline \Omega)$.
\label{thm2}
\end{thm}

Our proof of Proposition  \ref{thmmain}, 
 based on the general perturbation Lemma
3.1 in \cite{LN},   is similar to
 that of
Proposition 4.1 in
 \cite{LN}.

\noindent{\bf Proof of 
Proposition \ref{lem1}.}\
For any $\bar x\in B_1$, we see from (\ref{abab3}) that
as $r\to 0$,
\begin{eqnarray*}
\av _{ B_r(\bar x) }|u(x)- a(\bar x)|dx
&\le & 
\av  _{ B_r(\bar x) }|u(x)- [a(\bar x)+b(\bar x)
\cdot (x -\bar x)]|
dx
\\
&&+\av _{ B_r(\bar x) }
|b(\bar x)
\cdot (x -\bar x)|dx
\to 0.
\end{eqnarray*}
Thus, by a theorem of Lebesgue,
$
a=u$ a.e. in $ B_1.
$
We now take $
u\equiv a$, after changing the values of $u$ on
a zero measure set.
Let $\bar x, \bar y\in B_{1/4}$ satisfy, for some positive integer $k$,
$
2^{-(k+1)}\le |\bar x-\bar y|\le 2^{-k}.
$
By (\ref{abab3}), we have, for some dimensional constant $C$,
\begin{eqnarray*}
&&
|
u(\bar x) - [u(\bar y)+b(\bar y)
\cdot (\bar x -\bar y)]|
\\
&=&
|\av _{ B_{2^{-k}}(\bar x)  }
\{[u(\bar x)+b(\bar x)
\cdot (x -\bar x)]- [u(\bar y)+b(\bar y)
\cdot (x -\bar y)]\}
dx
|
\\
&\le&
\av
_{ B_{2^{-k}}(\bar x)  }
|u(x)- [u(\bar x)+b(\bar x)
\cdot (x -\bar x)]|
dx
\\
&&
+
\av _{ B_{2^{-k}}(\bar x)  }
|u(x)-  [u(\bar y)+b(\bar y)
\cdot (x -\bar y)]|
dx\\
&\le&\av 
_{ B_{2^{-k}}(\bar x)  }
|u(x)- [u(\bar x)+b(\bar x)
\cdot (x -\bar x)]|
dx\\
&& +2^n\av _{ B_{2^{-(k-1)}}(\bar y)  }
|u(x)-  [u(\bar y)+b(\bar y)
\cdot (x -\bar y)]|
dx\\
&\le& C 2^{ -k }
\delta(2^{-(k-1)})\le  C|\bar x-\bar y|\delta(4|\bar x-\bar y|).
\end{eqnarray*}
Switching the roles of $\bar x$ and $\bar y$ leads to
\begin{equation}
|
u(\bar y) - [u(\bar x)+b(\bar x)
\cdot (\bar y -\bar x)]|\le C|\bar y-\bar x|\delta(4|\bar y-\bar x|).
\label{ll1}
\end{equation}
Thus, by the above two inequalities and the triangle inequality,
\begin{equation}
|b(\bar x)-b(\bar y)|\le 2C \delta(4|\bar x-\bar y|).
\label{ll2}
\end{equation}
The conclusion of Proposition \ref{lem1} follows 
 from (\ref{ll1}) and (\ref{ll2}).
 
\medskip

\noindent{\bf Proof of Proposition \ref{thmmain}.}\
For simplicity, we prove it for $h=0, f=0$ --- the general case
only requires minor changes. 
We may assume
   without loss of generality that
$\varphi(1)\le \epsilon_0,
 \int_0^1 r^{-1}\varphi(r)dr\le \epsilon_0,
$
 for some small universal constant $\epsilon_0>0$.
This can be achieved by working with
$u(\delta_0x)$ for some $\delta_0$ satisfying
$\varphi(\delta_0)\le \epsilon_0$ and
$\int_0^{\delta_0}r^{-1}\varphi(r)dr<\epsilon_0$.
The smallness of $\epsilon_0$ will  be either obvious
or specified in the proof.
In the proof, a universal constant means
that it depends only on $\varphi$, $n, \lambda, \Lambda, N, \mu$.
We assume that $u$  is
normalized to satisfy
$
\|u\|_{  L^2(B_2) }=\varphi(4^{-1}).
$
We often write
$
\partial_\alpha(A_{ij}^{\alpha\beta}\partial_\beta u^j)
$
as
$
\partial(A\partial u).
$
For $k\ge 0$, let
$$
A_{k+1}(x)=A(4^{-(k+1)}x),\quad
\overline A=A(0).
$$
We will find $w_k\in H^1(\frac 3{4^{k+1}}
B_1, \Bbb R^N)$,  such that
for all $k\ge 0$,
\begin{equation}
\partial (\overline A\partial w_k) = 0, \   \frac 3{4^{k+1}} B_1,
\label{vk1}
\end{equation}
\begin{equation}
\|w_k\|_{  L^2(\frac 2{4^{k+1}}B_1) }\le C'4^{
- \frac{  k(n+2) }2   }\varphi(4^{-k}),\qquad
\|\nabla w_k\|_{  L^\infty(\frac 1{4^{k+1}}B_1) }\le C'
\varphi(4^{-k}),
\label{vk2}
\end{equation}
\begin{equation}
\|\nabla^2  w_k\|_{  L^\infty(\frac 1{4^{k+1}}B_1) }\le C'
4^k \varphi(4^{-k}),
\label{24new}
\end{equation}
\begin{equation}
\|u-\sum_{j=0}^{ k }w_j\|_{
L^2(  (\frac 14)^{k+1}  B_1)  }
\le
  4^{- \frac{(k+1)(n+2)}2 }
\varphi(4^{-(k+1)}).
\label{vk3}
\end{equation}
An easy consequence of (\ref{vk2}) is
\begin{equation}
\|w_k\|_{  L^\infty(4^{ -(k+1) }B_1)  }
\le  C'
4^{  -k }\varphi(4^{  -k }).
\label{qq1}
\end{equation}

\medskip

Here and in the following $C, C'$ and $\epsilon_0$ denote
various universal constants.
In particular,
they are independent of
$k$.
$C$ will be chosen first and will be large,
then
$C'$  (much larger than $C$),
and finally  $\epsilon_0\in (0,1)$
 (much smaller than $1/C'$).

 By Lemma 3.1 in \cite{LN},  we can find $w_0\in H^1(\frac 34
B_1, \Bbb R^N)$ such that
$$
\partial (\overline A\partial w_0) = 0, \
\qquad \mbox{in}\ \   \frac 34B_1,
$$
$$
\|u-w_0\|_{ L^2(\frac 12B_1) }\le C\epsilon_0^\gamma \|u\|_{ L^2(B_1)}
\le
4^{ -\frac{n+2}2 }\varphi(4^{-1}).
$$
So
$$\|w_0\|_{  L^2(\frac 12B_1) },\
\|\nabla w_0\|_{  L^\infty(\frac 14B_1) }, \
\|\nabla^2 w_0\|_{  L^\infty(\frac 14B_1) }\le C\varphi(1)\le C'
\varphi(1).
$$

We have verified (\ref{vk1})-(\ref{vk3}) for $k=0$.
Suppose that (\ref{vk1})-(\ref{vk3}) hold up to $k$ ($k\ge 0$);
 we will prove
them for $k+1$.
Let
$$
W(x)=  [ u-\sum_{j=0}^{ k }
w_j] (4^{ -(k+1) }x),
$$
$$
g_{k+1}(x)=4^{ -(k+1) }
\{  [\overline A-A_{k+1}](x)
\sum_{j=0}^k
 (\partial w_j)(4^{ -(k+1) }x)\}.
$$
Then $W$ satisfies
$$
\partial (A_{ k+1 }\partial W) =
\partial(g_{k+1}),
\qquad B_1.
$$
A simple calculation yields, using (\ref{BBB}),
\begin{eqnarray}
\|A_{k+1}-\overline A \|_{ L^2(B_1) }
&=& 
 \sqrt{|B_1|} \varphi(4^{-(k+1)})\le C(n,\mu) \varphi(4^{-(k+2)}).
\nonumber
\end{eqnarray}

By the induction hypothesis (see (\ref{vk2})-(\ref{vk3})),
$$
\sum_{j=0}^k
 |(\nabla w_j)(4^{ -(k+1) } x)|
\le C'\sum_{j=0}^k \varphi(4^{ -j})\le 
C(n) C' \int_0^1 r^{-1} \varphi(r)dr
 \le C(n)C'\epsilon_0,
\qquad x\in B_1,
$$
$$
\sum_{j=0}^k
 |(\nabla^2  w_j)(4^{ -(k+1) } x)|
\le 
 C'\sum_{j=0}^k 4^j \varphi(4^{ -j}),
\qquad x\in B_1,
$$
$$
\|W\|_{  L^2(B_1) }\le  4^{ -(k+1)  }\varphi(4^{ -(k+1)  })
\le C(\mu)  4^{ -(k+2) } \varphi(4^{ -(k+2)  }),
$$
$$
\|g_{k+1}\|_{ L^2(B_1) }
\le C(n,\mu)  C' \epsilon_0  4^{ -(k+2)  }\varphi(4^{ -(k+2)  }).
$$
By  Lemma 3.1 in \cite{LN},  there exists
$v_{k+1}\in H^1(\frac 34 B_1, \Bbb R^N)$ such that
$$
\partial(\overline A \partial v_{k+1})=0,
\qquad \ \mbox{in}\ \ \frac 34 B_1,
$$
and, for some universal constant $\gamma>0$, 
\begin{eqnarray}
\|W-v_{k+1}\|_{ L^2(\frac 12B_1) }&\le & C
(\|g_{k+1}\|_{ L^2(B_1) }+
 \epsilon_0^\gamma
\|W\|_{ L^2(B_1) })\nonumber\\
&\le& C   (C' \epsilon_0+\epsilon_0^\gamma)
4^{  -(k+2) } \varphi( 4^{  -(k+2) } ).
\label{vk4}
\end{eqnarray}

Let
$$
w_{k+1}(x)=v_{k+1}(4^{k+1}x), \qquad x\in \frac{ 3 }{  4^{k+2} }B_1.
$$
A change of variables in (\ref{vk4}) and in the equation of
$v_{k+1}$  yields
 (\ref{vk1}) and (\ref{vk3}) for $k+1$.
It follows from the above that
$$
\|\nabla ^2  v_{k+1}\|_{ L^\infty(\frac 14B_1) }+
\|\nabla v_{k+1}\|_{ L^\infty(\frac 14B_1) }\le  C
\|v_{k+1}\|_{ L^2(\frac 12B_1) }\le  C
4^{  -(k+1) } \varphi( 4^{  -(k+1) } ).
$$
Estimates
(\ref{vk2}) for $k+1$ follow from the above estimates
for $v_{k+1}$.  We have, thus,   established
(\ref{vk1})-(\ref{vk3}) for all $k$.

\medskip

For $x\in   4^{ -(k+1)}B_1$,
  using (\ref{vk2}), (\ref{24new}),  (\ref{qq1}), (\ref{BBB}) and
 Taylor expansion,
\begin{eqnarray}
&&| \sum_{j=0}^k
w_j(x)- \sum_{j=0}^\infty
w_j(0)  -\sum_{j=0}^\infty
\nabla w_j(0)\cdot x|\nonumber\\
&\le &
 \sum_{j=k+1}^\infty(|w_j(0)|+|\nabla w_j(0)||x|)
+ \sum_{j=0}^k
\|\nabla^2 w_j\|_{ L^\infty(4^{-(k+1) }B_1)}|x|^2\nonumber\\
&\le &
C  \sum_{j=k+1}^\infty  (4^{-j}\varphi(4^{-j})
+\varphi(4^{-j})|x|)+C \sum_{j=0}^k
4^j\varphi(4^{-j})|x|^2\nonumber\\
&\le &
C  4^{ -(k+1)  }   
  \int_0^{4^{-k}} r^{-1}  \varphi(r)dr+
C |x|^2 \int_{ \frac {|x|}2 }^1 r^{-2}  \varphi(r)
dr.  
\end{eqnarray}
It is easy to see 
that $\lim_{|x|\to 0} |x|\int_{ \frac {|x|}2 }^1 r^{-2}  \varphi(r)
dr=0$, since (\ref{BBB}) implies $\lim_{r\to 0^+}\varphi(r)=0$.

We then derive from (\ref{vk3}) and the above, using H\"older inequality,
  that, for some $\delta(r)=\circ(1)$ (as $r\to 0$),  
 depending only on $
\varphi$, $n, \lambda, \Lambda, N, \mu$,
\begin{eqnarray}
&& \int_{ 4^{ -(k+1)}  B_1 }|u(x)- ( \sum_{j=0}^\infty
w_j(0) +\sum_{j=0}^\infty
\nabla w_j(0)\cdot x)|dx\nonumber
\\
&\le & 
 \| \sum_{j=0}^k
w_j(x)- \sum_{j=0}^\infty
(w_j(0)  -
\nabla w_j(0)\cdot x)\|_{ L^1(4^{ -(k+1)}B_1) }
+ \|u-  \sum_{j=0}^k
w_j(x)\|_{ L^1(4^{ -(k+1)}B_1) }\nonumber\\ 
&=&  4^{-(k+1)(n+1)} \delta(4^{-(k+1)}).
\nonumber
\end{eqnarray}
Proposition \ref{thmmain} follows from the above with
$a=\sum_{j=0}^\infty
w_j(0)$ and $b=\sum_{j=0}^\infty
\nabla w_j(0)\cdot x)$.

\noindent{\bf Proof of Theorem \ref{thm1}.}\
Fix a $\rho\in C^\infty_c(B_4)$, $\rho\equiv 1$ on $B_3$, and let
$$
\varphi(r):=
\sup_{x\in B_3}
( \av_{B_r(x)} |(\rho A)- (\rho A)(x)|^2)^{\frac 12}.
$$
It is easy to see that for some $\mu>1$, $\varphi$ satisfies
(\ref{BBB}). 
Indeed, since it is easily seen that 
$$
\varphi(r)\le C(\bar \varphi(r)+r),
$$
the second inequality follows.
For the first inequality,  
we only need to show that $\varphi(2r)\le C(n)\varphi(r)$, since the rest 
is obvious.
For any $\bar x$, let $x_1=\bar x, x_2,  \cdots, x_m$, $m=m(n)$, satisfy
$B_{2r}(\bar x)\subset \cup_{i=1}^m B_{r/9}(x_i)$, and
$|x_i-x_{i+1}|\le r/9$.
Then
\begin{eqnarray*}
&&( \av_{B_{2r}(\bar x)} |(\rho A)-(\rho A)(\bar x)|^2)^{\frac 12}
\\&\le& C(n)\sum_{i=1}^m
( \av_{B_{r/9}(x_i)} |(\rho A)-(\rho A)(\bar x)|^2)^{\frac 12}\\
&\le& C(n) \sum_{i=1}^m
\{( \av_{B_{r/9}(x_i)} |(\rho A)-(\rho A)(x_i)|^2)^{\frac 12}
+
|(\rho A)(\bar x)-(\rho A)(x_i)|\}\\
&\le &  C(n) \varphi(r)
+C(n)  \sum_{i=1}^{m-1}
 |(\rho A)(x_i)-(\rho A)(x_{i+1})|.
\end{eqnarray*}
Since
\begin{eqnarray*}
&& |(\rho A)(x_i)-(\rho A)(x_{i+1})|\\
&=& | \av_{B_{r/9}(x_i)} [(\rho A)-(\rho A)(x_i)]
-  \av_{B_{r/9}(x_i)} [(\rho A)-(\rho A)(x_{i+1})]|
\\
&\le & C(n) ( \av_{B_{r}(x_i)} |(\rho A)-(\rho A)(x_i)|
+  \av_{B_{r}(x_{i+1})}  |(\rho A)-(\rho A)(x_{i+1})|)\le
C(n)\varphi(r),
\end{eqnarray*}
we have
$$
( \av_{B_{2r}(\bar x)} |(\rho A)-(\rho A)(\bar x)|^2)^{\frac 12}
\le C(n)\varphi(r).
$$
Thus  $\varphi(2r)\le C(n)\varphi(r)$.

For any $\bar x\in B_2$, 
$$
( \av_{B_r(\bar x)} |A-A(\bar x)|^2)^{\frac 12}\le
\varphi(r),\qquad 0<r<1/4.
$$

Thus Theorem \ref{thm1} follows from 
Proposition \ref{thmmain}-\ref{lem1}.

\end{document}